\documentclass{amsart}
\usepackage{amssymb,latexsym}
\numberwithin{equation}{section}

\newtheorem{theorem}{Theorem}[section]
\newtheorem{lemma}[theorem]{Lemma}
\newtheorem{proposition}[theorem]{Proposition}
\newtheorem{corollary}[theorem]{Corollary}

\theoremstyle{definition}

\theoremstyle{remark}

\def\N{{\mathbb N}}
\def\Z{{\mathbb Z}}%
\def\F{{\mathbb F}}
\def\Q{{\mathbb Q}}
\def\C{{\mathbb C}}%

\def\dim{{dim}}\def\ker{{Ker}}\def\udim{\underline{{\dim}}}
\def\hom{{Hom}}\def\im{{Im}}\def\endo{{End}}
\def\ext{{Ext^1}}
\def\aut{{\hbox{Aut}}}
\def\irj{i\rightarrow j}
\def\bd{{\mathbf d}}
\def\qo{Q^{op}}\def\po{p^{op}}
\def\vd{V_d}\def\rd{R_d}\def\gd{G_d}
\def\fdb{{\mathcal F}_{\mathbf d}}\def\pdb{{P}_{\mathbf d}}
\def\udb{{U}_{\mathbf d}}\def\ldb{{L}_{\mathbf d}}
\def\xdb{{X}_{\mathbf d}}\def\ydb{{Y}_{\mathbf d}}
\def\txdb{{\tilde X}_{\mathbf d}}

\def\ox{{\mathcal O}_X}\def\om{{\mathcal O}_M}\def\on{{\mathcal O}_N}
\def\xn{X_{\overline N}}\def\xnm{X_{\overline N}^{\overline M}}
\def\txn{\tilde X_{\overline N}}\def\txnm{\tilde X_{\overline N}^{\overline M}}
\def\yon{Y_{\overline N}}\def\yn{Y_{N}}
\def\lgd{{\mathfrak g}_d}\def\ludb{{\mathfrak u}_{{\mathbf d}}}
\def\hp{{F}_{N_\nu,\ldots,N_1}^M}\def\hvq{{\mathcal H}_v(Q)}

\begin{document}
\title[The bar automorphism and the preprojective variety]
{The bar automorphism in quantum groups and geometry of quiver representations}
\author{Philippe Caldero}
\address{D\'epartement de Math\'ematiques, Universit\'e Claude Bernard Lyon I,
69622 Villeurbanne Cedex, France}
\email{caldero@igd.univ-lyon1.fr}
\author{Markus Reineke}
\address{Mathematisches Institut, Universit\"at M\"unster, 48149 M\"unster, 
Germany}
\email{reinekem@math.uni-muenster.de}

\begin{abstract}
Two geometric interpretations of the bar automorphism in the positive part of a 
quantized enveloping algebra are given.
The first is in terms of numbers of rational points over finite fields
of quiver analogues of orbital varieties; the second is in terms of a duality
of constructible functions provided by preprojective varieties of quivers.
\end{abstract}

\date{\today}

\maketitle
\begin{section} {Introduction}
The canonical basis $\mathcal{B}$ of the positive part $U_v(\mathfrak g)^+$ of 
the quantized enveloping algebra of a semisimple Lie algebra $\mathfrak{g}$, 
constructed by G.~Lusztig \cite{lusztig1}, has many favourable properties, like 
for example inducing bases in all the finite dimensional irreducible 
representations of $\mathfrak{g}$ simultaneously.\par
The basis $\mathcal{B}$ can be characterized algebraically by its elements being 
fixed under the so-called bar automorphism of $U_v(\mathfrak g)^+$, and by 
admitting a unitriangular base change to the PBW-type bases.\par
The Hall algebra approach to quantum groups \cite{ringelinv} provides a 
realization of certain specializations $U_q(\mathfrak g)^+$ via a convolution 
product for constructible functions on varieties $\rd$ parametrizing 
representations of Dynkin quivers. In this realization, the elements of PBW-type 
bases correspond to characteristic functions of orbits in $\rd$ under a natural 
action of an algebraic group $\gd$, whereas the elements of $\mathcal B$ 
correspond to constructible functions naturally associated to the intersection 
cohomology complexes of the closures of $\gd$-orbits.\par 
It is therefore natural to also ask for interpretations of the bar automorphism 
in terms of the geometry of the varieties $\rd$, since this automorphism plays a 
central role in defining the canonical basis algebraically.\par
In the present paper, two such interpretations are given. In the geometric setup 
of \cite{lusztigbook}, analogues of orbital varieties in the quiver context are 
constructed. These parametrize quiver representations fixing certain flags, and 
their numbers of rational points over finite fields are shown (Theorem 
\ref{main}) to give essentially the coefficients $\Omega_{M,N}$ of the bar 
automorphism on a PBW-type basis of $U_q(\mathfrak g)^+$. The key ingredient in 
deriving this result in section 3 is a generalization (Corollary 
\ref{filtration}) of a very useful formula \cite{riedtmann1} by C.~Riedtmann, 
relating numbers of filtrations of quiver representations over finite fields to 
cardinalities of orbital varieties; this generalization, together with the 
construction of the orbital varieties, is given in section 2.\par
The second interpretation starts from a duality between constructible functions 
on the varieties $\rd$ provided by the preprojective varieties of a quiver, also 
used by G. Lusztig \cite{lusztig2}. The coefficients $\Omega_{M,N}$ are shown 
(Proposition \ref{tdual}) to be essentially given by a convolution operator 
derived from a certain twisted version of this duality, constructed in section 
4.\par
{\bf Acknowledgments:} The present work was done during a stay of both authors 
at the University of Antwerp in the framework of the Priority Programme "Non-
Commutative Geometry" of the European Science Foundation.

\end{section}
\begin{section} {A generalization of Riedtmann's formula}
\begin{subsection}{ }

Let $k$ be a field. We fix a finite quiver $Q$ with 
set 
of vertices $I$ and set of arrows $Q_1$, whose underlying unoriented graph is a 
disjoint union of Dynkin diagrams of type $A_n$, $D_n$, $E_6$, $E_7$ or $E_8$. 
Fixing a dimension type $d=\sum_{i}d_ii\in\N I$,  define an 
$I$-graded vector space $\vd:=\oplus_{i\in I}k^{d_i}$.
For any subquotient $W$ of $\vd$ compatible with the grading, we set
$\udim(W):=\sum_i(\dim W_i)i\in\N I$, where $W_i$ denotes the $i$-component of $W$.
\par
Set $\rd:=\rd(Q):=\bigoplus_{\alpha:\irj}\hom(k^{d_i},k^{d_j})$ and 
$\gd:=\prod_i GL_{d_i}(k)$. 
The algebraic group $\gd$ acts linearly on the affine space $\rd$ by 
$(g.M)_{\alpha:\irj}:=g_jM_\alpha g_i^{-1}$. 
   
\end{subsection}
\begin{subsection}{ }

Let $\nu$ be a positive integer, and let $\bd=(d^1,\ldots,d^s)$ be a $\nu$-tuple 
of elements 
in $\N I$ such that $d=\sum_s d^s$. Let $\fdb$ be the set of filtrations $F^*=
(0=F^\nu\subset\ldots\subset F^1\subset F^0=\vd)$
of the graded space $\vd$ such that $\udim F^{s-1}/F^s=d^s$, $1\leq s\leq\nu$.
The action of the group $\gd$ on $\vd$ provides a transitive action of $\gd$ on  
$\fdb$.
\par
We fix an arbitrary filtration $F_0^*$ in $\fdb$. Choosing successive 
complements, we can assume that $\vd$ has a direct sum decomposition 
$\vd=\bigoplus_{s=1}^\nu\vd^s$ (as $I$-graded $k$-space), such that 
$F_0^s=\bigoplus_{t>s}\vd^t$ for $s=1,\ldots,\nu$. This induces a decomposition 
$\rd=\bigoplus_{s,t=1}^\nu\rd^{s,t}$ by setting 
$\rd^{s,t}=\bigoplus_{\alpha:\irj}\hom((\vd^s)_i,(\vd^t)_j)$. Let $\pdb$ be the 
stabilizer of $F_0^*$ in $\gd$. This
is a 
parabolic subgroup of $\gd$, providing an identification $\fdb\simeq\gd/\pdb$.
\par
We say that an element $M$ in $\rd$ is compatible with a filtration $F^*$ in 
$\fdb$ if $M_{\alpha}
(F_i^s)\subset F_j^s$ for any arrow $\alpha:\irj$ in $Q_1$ and any 
$s=1,\ldots,\nu$. Set 
$$\xdb:=\{(M,F^*),\,M\in\rd,\,F^*\in\fdb,\,
\hbox{$M$ compatible with $F^*$}\}\subset \rd\times\fdb.$$
The diagonal action of the group $\gd$ on $\rd\times\fdb$ respects $\xdb$, and
the projections $p_1$ and $p_2$ on $\rd$ and $\fdb$, respectively, are 
$\gd$-equivariant. \par
Set $\ydb:=p_2^{-1}(F_0^*)$, which will be identified with its image 
$p_1(\ydb)$ in $\rd$. We have the identification 
\begin{equation}\label{identification}
 \gd\times^{\pdb}\ydb\simeq\xdb,\,\,\overline{(g,y)}:=\{(gp,p^{-1}.y),\,
p\in\pdb\}\mapsto (g.y,gF_0^*).
\end{equation}
Denote by $\udb$ the unipotent radical of $\pdb$ and set $\txdb:=
\gd\times^{\udb}\ydb$. Let $\pi$ be the natural projection $\txdb\rightarrow
\gd\times^{\pdb}\ydb\simeq\xdb$. The Levi decomposition gives
$\pdb=\ldb\udb$, where $\ldb$ is the Levi of the parabolic $\pdb$.
We have a diagonal action of $\ldb$ on $\txdb$ defined by $l.(g,y):=
(gl^{-1},l.y)$, and a left action of $\gd$ on $\txdb$. These 
two actions commute.\par
The above direct sum decomposition $\vd=\bigoplus_s\vd^s$ provides a surjection 
$\zeta:\rd\rightarrow
\prod_s R_{d^s}$, which restricts to a surjection $\zeta$ : 
$\ydb\rightarrow\prod_s R_{d^s}$. This clearly defines a surjection 
$\tilde \zeta$ : 
$\gd\times^{\udb}\ydb\rightarrow\prod_s R_{d^s}$ by projecting on the second 
factor.
We obtain the following diagram.
\begin{equation}
\begin{matrix}\prod_s R_{d^s}&{\stackrel{\tilde 
\zeta}{\leftarrow}}&\txdb&{}&{}\cr
{}&{}&\downarrow\pi&&{}&{}\cr
\rd&{\stackrel{p_1}{\leftarrow}}&\xdb&{\stackrel{p_2}{\rightarrow}}&\fdb\simeq
\gd/\pdb\cr
\end{matrix}
\end{equation}

The following lemma follows immediately from the definitions:
\begin{lemma}\label{commute} 
The morphism $\tilde \zeta$ commutes with the action of $\ldb$.
 The morphisms $\pi$, $p_1$ and $p_2$ commute with the action of $\gd$.
\end{lemma}
\end{subsection}
\begin{subsection}{}

Let $\bmod kQ$ be the category of finite dimensional $k$-representations of $Q$. 
For a representation $X$ in 
$\bmod kQ$ of dimension type $\udim(X)=d$, we denote by $\ox$ the corresponding 
$\gd$-orbit in $\rd$. The $\gd$-orbits are in bijection with the isoclasses of 
representations of dimension type $d$ in $\bmod kQ$ by definition. We  denote by 
$\overline{\bmod }kQ$ this set of isoclasses, and $\overline X$ will be the 
isoclass
corresponding to the representation $X$.\par
We now suppose that $\nu$ is the number of isoclasses of indecomposable 
representations in 
$\bmod kQ$ (which coincides with the number of positive roots of the root system 
corresponding to $Q$ by Gabriel's Theorem), and we denote by $I_s$ for $1\leq 
s\leq\nu$ these 
indecomposable representations, ordered such that 
\begin{equation}\label{directed}
\hom_Q(I_t,I_s)=0\mbox{ for } 1\leq s<t\leq\nu,
\end{equation}
where $\hom_Q$ denotes the space of homomorphisms in the category $\bmod kQ$ 
(such an ordering exists since the category $\bmod kQ$ is directed; for this and 
other facts on $\bmod kQ$ see, for example, \cite{ringelbook}).\par
We fix two representations $N$ and $M$ in $\rd$. In the following, 
we define analogues of
 the varieties introduced in the previous section, naturally associated to $N$ 
and 
$M$.\par
Let $N=\oplus_s N_s$ be the unique decomposition of the representation $N$ 
such that $N_s$ is isomorphic to a direct sum of copies of $I_s$ for $1\leq 
s\leq\nu$.
We can 
suppose without loss of generality that the spaces $N_s$ are compatible with 
the direct sum decomposition $\vd=\bigoplus_s\vd^s$, so that $N_s$ belongs to 
$\rd^{s,s}$.\par
We define  $\xn$ as the set of pairs $(P,F^*)$ in $\xdb$ such that the 
representation induced by $P$ on $F^{s-1}/F^s$ is isomorphic to $N_s$ for any 
$s=1,\ldots,\nu$. 
Let $\xnm$ be the subset of pairs $(P,F^*)$ in $\xn$ such that $P$ belongs to 
$\om$.
Set $\yon:=\zeta^{-1}(\prod_s{\mathcal O}_{N_s})\subset\ydb$. Again, the 
following 
lemma follows immediately from the definitions:
\begin{lemma}\label{xn}
Via the identification \ref{identification}, we have
$$\xn\simeq\gd\times^{\pdb}\yon\mbox{ and }\xnm\simeq\gd\times^{\pdb}
(\yon\cap\om).$$ 

\end{lemma}

We define
$$\txnm:=\gd\times^{\udb}(\yon\cap\om)\subset\txn:=\gd\times^{\udb}\yon
\subset
\xdb.$$
The left action of $\gd$ and the diagonal action of $\ldb$ both stabilize these 
varieties.
\end{subsection}
\begin{subsection}{}

From now on, we suppose that $k$ is a finite field with $q$ elements.
For any representation $P$ in $\rd$, we denote by $\aut(P)\subset\gd$ the 
stabilizer of $P$ (which coincides with the automorphisms of $P$ as an object in 
$\bmod kQ$ by definition). With the notation of the previous section, 
let $\yn$ be the fiber of $\zeta$ over the point $(N_s)$, that is, $\yn:=
\zeta^{-1}((N_s))\subset\ydb$. Note that $\udb$ acts on $\yn$.
\begin{proposition}\label{fiber} 
With notation as above, we have
$$|p_1^{-1}(M)\cap\xn|=\frac{|\aut(M)|\cdot|\yn\cap\om|}
{(\prod_s|\aut(N_s)|)\cdot|\udb|}.$$
\end{proposition}
\begin{proof}
First we have 
 $$ \pi^{-1}p_1^{-1}(\om)\cap\txn=\txnm\simeq\gd\times^{\udb}(\yon\cap\om).$$
Since $p_1\pi$ commutes with the $\gd$-action, we obtain
\begin{eqnarray*}|\pi^{-1}p_1^{-1}(M)\cap\txn|&=&\frac{|\gd|\cdot|\yon\cap\om|}
{|\om|\cdot|\udb|}\\
&=&\frac{|\aut(M)|\cdot|\yon\cap\om|}{|\udb|}.
\end{eqnarray*}
From Lemma \ref{commute}, we conclude 
\begin{eqnarray*}
|p_1^{-1}(M)\cap\xn|&=&\frac{|\aut(M)|\cdot|\yon\cap\om|}{|\ldb|\cdot|\udb|}\\
&=&\frac{|\aut(M)|\cdot(\prod_s|{\mathcal O}_{N_s}|)\cdot|\yn\cap\om|}{|\ldb
|\cdot|\udb|}.\end{eqnarray*}
This implies the proposition.
\end{proof}
\end{subsection}
\begin{subsection}{}

In this section, we give a more precise version of Proposition \ref{fiber}. Let
$\lgd:=\bigoplus_i\mathfrak{gl}(k^{d_i})$ be the Lie algebra of the group $\gd$.
The components of 
an element $\xi$ in $\lgd$ will be denoted by $\xi_i$. Let $\ludb\subset\lgd$ 
be the Lie algebra of $\udb$. The differential of the morphism 
$\gd\rightarrow\gd.N$ 
gives rise to a morphism of vector spaces $\phi$ : 
$\lgd\rightarrow\rd$ given by 
$\phi(\xi)_{\alpha:\irj}=\xi_jN_{\alpha}- N_{\alpha}\xi_i$. 

\begin{lemma}\label{phi}
The morphism $\phi$ has the following properties :
\item{(i)} $\ker(\phi)=\endo_Q(N)$,\par
\item{(ii)} $\im(\phi)=T_N$, where $T_N:=T_N(\on)$ is the tangent space to 
$\on$ at the 
point $N$,\par
\item{(iii)} $\im(\phi)$ is compatible with the decomposition 
$\rd=\bigoplus_{s,t}\rd^{s,t}$, and it contains the subspace 
$\bigoplus_{s\leq t}\rd^{s,t}$,\par
\item{(iv)} the restriction of $\phi$ to $\ludb$ is injective.
\end{lemma}
\begin{proof}
(i) follows from the definition of $\phi$ and of the category $\bmod kQ$. 
(ii) is clear. The first assertion
of (iii) follows from the fact that $\phi$ decomposes into a direct sum 
$\phi=\bigoplus_{s,t}\phi_{s,t}$, where
\begin{equation}\label{decompphi}
\phi_{s,t}\,:\,\bigoplus_i\hom_k((\vd^s)_i,(\vd^t)_i)
\rightarrow\rd^{s,t}.
\end{equation}
The second assertion is \cite[Lemma 10.4]{lusztig1}.
To prove (iv), we remark that $\phi|_{\ludb}=
\bigoplus_{s>t}\phi_{s,t}$ which implies that $\ker\phi|_{\ludb}=0$ by 
formula \ref{directed} and (i).
\end{proof}
We now consider the affine space $\yn$ identified with its tangent vector 
space $T_N(Y_N)$ by $x\mapsto x+N$.

\begin{proposition}\label{slice}
The space $\yn$ contains $\phi(\ludb)$. Let $E_N$ be any complement of 
$\phi(\ludb)$ in $\yn$
which is compatible with the decomposition 
$\bigoplus\rd^{s,t}$ of $\rd$. Then,\par
\item{(i)} the space $E_N$ is a direct summand of $T_N(\on)$ in $\rd$,\par
\item{(ii)} $E_N$ is a tranversal slice for the action of $\udb$ on $\yn$.
\end{proposition}

\begin{proof}
First we remark that, by construction, the space $\yn$ is compatible with the 
decomposition $\rd=\bigoplus_{s,t}\rd^{s,t}$. Thus, a complement $E_N$ as above 
exists by Lemma 
\ref{phi} (iii).\par
(i) is a consequence of Lemma \ref{phi} (ii), (iii) and the decomposition 
\ref{decompphi} of $\phi$.\par
Now we prove (ii). Fix a complement $E_N$ and let $X$ be in $E_N$. We claim that 
$X$ is the 
unique element of $E_N$ in the orbit $\udb.X$. Suppose 
that $Y=U.X\in E_N$, with $U\in\udb$. The component of $U$ (induced by the 
decomposition $\vd=\bigoplus_s\vd^s$) belonging to 
$\bigoplus_i\hom((\vd^s)_i,(\vd^t)_i)$ is denoted by $U_{s,t}$ for $s\geq t$. 
Note that $U_{s,s}$ is the identity for all $s=1,\ldots,\nu$.\par
We prove by induction on 
$s-t>0$ that $U_{s,t}=0$, for $s>t$. Fix a pair $(s,t)$ for $1\leq t<s\leq\nu$.
It is easily seen by a weight argument that the induction hypothesis implies 
that the component $Y_{s,t}$ of $Y$ has the following form :
$Y_{s,t}=U_{s,t}N_t -N_s U_{s,t}$. So, we have $U_{s,t}N_t -N_s U_{s,t}=Y_{s,t}
\in E_N$ by the hypothesis on $E_N$. This implies that $Y_{s,t}\in E_N\cap
\phi(\ludb)=\{0\}$. By Lemma \ref{phi}, this gives $U_{s,t}=0$. The claim is 
proved.\par
In particular, this implies that the action of the group $\udb$ on $\yn$ is 
free.
Thus, $|E_N|=\frac{|\yn|}{|\udb|}=|\yn/\udb|$. This equality, together with the 
claim just proved, gives (ii).
\end{proof}
The following corollary can be seen as a generalization of Riedtmann's 
formula, \cite{riedtmann1}.
\begin{corollary}\label{filtration}
Let $E_N$ be as in the previous proposition. Then,
$$|{\mathcal F}_{N_\nu,\ldots,N_1}^M|=\frac{|\aut(M)|\cdot|E_N\cap\om|}
{\prod_s|\aut(N_s)|},$$
where ${\mathcal F}_{N_\nu,\ldots,N_1}^M$ denotes the set of filtrations 
$0=M^\nu\subset\ldots\subset M^1\subset M^0=M$ of 
the representation $M$ with successive subquotients $M^{s-1}/M^s$ isomorphic to 
$N_s$ for  
$1\leq s\leq\nu$.
\end{corollary}

\begin{proof}
By construction, we have ${\mathcal F}_{N_\nu,\ldots,N_1}^M=p_1^{-1}(M)\cap\xn$. 
Moreover, as $\udb\subset\gd$, the previous proposition implies that 
$E_N\cap\om$ is a tranversal slice for the action of $\udb$ on $\yn\cap\om$.
The corollary then follows from Proposition \ref{fiber}.
\end{proof}

It is known \cite{ringelhall} that there exists a polynomial $\hp(t)\in\Z[t]$, 
called the generalized Hall polynomial, whose value at any prime power $q$ 
equals
the number of ${\mathbb F}_q$-rational points of the variety ${\mathcal 
F}_{N_\nu,\ldots,N_1}^M$.
\end{subsection}
\end{section}

\begin{section}{Hall algebras and coefficients of the bar automorphism.}
\begin{subsection}{}

We define the Euler form $<\,,\,>$ on $\N I$ by $<d,e>:=\sum_{i\in I}d_ie_i-
\sum_{\alpha:\irj}d_ie_j$. We suppose in this section that $k$ is a field with 
$q=v^2$ 
elements for some $v\in\C$. The dimension type $d$ and the representations $M$, 
$N$
are no longer fixed. For all finite sets $X$ on which $\gd$ acts, we denote by 
$\C_{\gd}[X]$ the set of $\gd$-invariant functions from $X$ to $\C$.
Define
$$\hvq=\bigoplus_{d\in\N I}\C_{\gd}[\rd].$$
The space $\hvq$ is endowed with a structure of $\N I$-graded $\C$-algebra by 
the convolution product :
$$(f.g)(X)=v^{<d,e>}\sum_{U\subset X}f(X/U)g(U),\,\,f\in\C_{\gd}[\rd],\,
g\in\C_{G_e}[R_e],\,X\in R_{d+e},$$
where $U$ runs over all subrepresentations of $X$ of dimension type $e$.
It is known \cite{ringelhall} that this product defines the structure of an 
associative algebra on 
$\hvq$, which is called the (twisted) Hall algebra of the quiver $Q$.\par
For any representation $M$ of $Q$ with isoclass $\overline M$, let 
$e_M=e_{\overline M}$ be $v^{\dim\endo N-\dim N}$ times the 
characteristic 
function of the orbit $\om$. It is clear that 
$\{e_{\overline M}\}_{\overline M\in\overline{\bmod }kQ}$ is a basis of $\hvq$.
\par
Let $S_i$ be the simple representation corresponding to the vertex $i$ in $I$.
It is known \cite{ringelinv} that there exists an isomorphism $\eta$ from the 
Hall algebra 
$\hvq$ to the positive part $U_q({\mathfrak g})^+$ of the quantum enveloping 
algebra associated to $Q$, such that $\eta$ maps $e_{S_i}$ to the canonical 
generator $e_i$ of  $U_q({\mathfrak g})^+$. Note that the basis 
$e_M=e_{\overline M}$ is sent to the so-called Poincar\'e-Birkhoff-Witt basis
of $U_q({\mathfrak g})^+$ which corresponds to a reduced decomposition of the 
longest Weyl group element naturally associated to $Q$, see 
\cite[4.12]{lusztig1}.\par
\end{subsection}

\begin{subsection}{}

We consider the inner product on $\hvq$, called Green form, defined by
$$(e_M,e_N)=v^{2\dim\endo N}a_M^{-1}\delta_{N,M},$$
where $\delta$ is the Kronecker symbol and 
$a_M:=a_M(v^2)=|\aut(M)|$. Hence, we 
obtain the dual PBW type basis by setting $e_M^*=v^{-2\dim\endo N}a_Me_M$.

The following lemma is an easy consequence of the 
definition of the convolution product in the Hall algebra and of the properties 
of the decomposition $N=\bigoplus_sN_s$.

\begin{lemma}\label{decomposition}
Suppose that $N=\bigoplus_s N_s$ is the decomposition of $N$ into 
powers of indecomposables as above. Then,
\par
\item{(i)} $e_N=e_{N_1}\ldots e_{N_\nu}$,\par
\item{(ii)} $e_N^*=e_{N_1}^*\ldots e_{N_\nu}^*$,\par
\item{(iii)} $e_{N_\nu}\ldots e_{N_1}=\sum_M v^{S-\dim\endo M}\hp(v^2)e_M$,
where $S=\sum_s\dim\endo N_s+\sum_{s>t}<\udim N_s,\udim N_t>$.
\end{lemma}
\end{subsection}
\begin{subsection}{}

Using the basis elements $e_M$ introduced above, the multiplication in the Hall 
algebra reads as follows:
$$e_M\cdot e_N=\sum_{\overline X}v^{\dim\endo M+\dim\endo N+<\udim M,\udim N>-
\dim\endo X}F_{M,N}^X(v^2)\cdot e_X.$$
Since the $F_{M,N}^X$ are polynomials, we can thus take the above formula as the 
definition of structure constants for a $\Q(v)$-algebra, the generic (twisted) 
Hall algebra \cite{ringelinv}, which will still be denoted by $\hvq$.\par
We define a $\Q$-linear involution on $\Q(v)$ by $\overline v=v^{-1}$. We define 
on $\hvq$ :\par 
-- a $\overline{\cdot}$-linear involution by $\overline{e_i}=e_i$, called the 
bar involution,\par
-- the $\Q(v)$-linear antiinvolution $\sigma$ by $\sigma(e_i)=e_i$.\par
Denote by $\omega_{M,N}$ the $e_N$-coefficient of  $\overline{e_M}$ in the 
PBW-basis. It is clear that $\omega_{M,N}$ is zero if $M$ and $N$ do not have 
the same dimension type. Following Lusztig \cite{lusztig1}, 
we use the normalization $\Omega_{M,N}=v^{\dim\on -\dim\om} \omega_{M,N}
=v^{\dim\endo M -\dim\endo N}\omega_{M,N}$. We say that $M$ degenerates to $N$ 
if $\on$ belongs to the closure of $\om$. The following is proved in 
\cite{lusztig1}.

\begin{lemma}\label{omega}
For any two representations $M$, $N$ in $\bmod kQ$, we have:\par
\item{(i)} if $\Omega_{M,N}\not= 0$, then $M$ degenerates to $N$,\par
\item{(ii)} $\Omega_{M,M}=1$,\par
\item{(iii)} $\Omega_{M,N}\in\Z[v^{-2}]$.
\end{lemma}
We now want to give a geometric interpretation of the polynomial 
$\overline{\Omega_{M,N}}$. This will be provided by Theorem \ref{main}.
The following proposition precises a result of 
\cite[Proposition 3.1]{caldschi1}, where the formula was asserted up 
to a power of $v$.

\begin{proposition}\label{formula}
Let $M$, $N$ be two representations in $\bmod kQ$, and let 
$N=\bigoplus_s N_s$ be the decomposition into powers of indecomposables as 
above. Then, the polynomial $\overline{\Omega_{M,N}}\in\Z[v^2]$ 
is given by 
$$\overline{\Omega_{M,N}}=\hp(v^2)\frac{\prod_s a_{N_s}}{a_M}.$$
\end{proposition}

\begin{proof}
We first calculate $\omega_{M,N}=(\overline{e_M},e_N^*)$ by using the adjoint 
of the bar automorphism for the Green form.
From \cite[1.2.10.]{lusztigbook}, we have :
$$\overline{\omega_{M,N}}=(-v)^{-\dim M}v^{-<\udim N, 
\udim N>}(e_M,\sigma(\overline{e_N^*})).$$
From Lemma \ref{decomposition} (ii), this gives :
 $$\overline{\omega_{M,N}}=(-v)^{-\dim M}v^{-<\udim N, 
\udim N>}(e_M,\sigma(\overline{e_{N_\nu}^*})\ldots
\sigma(\overline{e_{N_1}^*})).$$
By \cite{caldero1}, the elements $e_{N_s}^*$ belong to the dual canonical basis.
Hence, by \cite[Lemma 4.3.]{reineke1},
$$\sigma(e_{N_s}^*)=(-v)^{\sum\dim N_s}v^{<\udim N_s,\udim N_s>}e_{N_s}^*$$
for all $s=1,\ldots,\nu$. We deduce that 
$$\overline{\omega_{M,N}}=v^{-\sum_{s\not=t}<\udim N_s, 
\udim N_t>}(e_M,e_{N_\nu}^*\ldots e_{N_1}^*).$$
It remains to calculate the $e_M^*$-component of the product 
$e_{N_\nu}^*\ldots e_{N_1}^*$ in the dual PBW-basis. This is obtained from
Lemma \ref{decomposition} :
$$(e_M,e_{N_\nu}^*\ldots e_{N_1}^*)=v^T\hp a_{N_1}\ldots a_{N_\nu}a_M^{-1},$$
where $T=-\sum_s\dim\endo N_s+\dim\endo M+\sum_{s>t}<\udim N_s,\udim N_t>.$
Now, from the interpretation of $<\,,\,>$ as the homological Euler form in 
$\bmod kQ$, namely $<\udim M,\udim N>=\dim\hom(M,N)-\dim\ext(M,N)$, we easily 
obtain
$$\dim\endo(N)=\sum_s\dim\endo(N_s)+\sum_{s\leq t}<\udim N_s,\udim N_t>,$$
and the claimed formula follows.

\end{proof}
\end{subsection}
\begin{subsection}{}

Our efforts are rewarded in that we can deduce a geometric interpretation 
of the coefficient $\overline{\Omega_{M,N}}$.

\begin{theorem}\label{main} 
Let $k$ be the finite field $\F_{v^2}$. Fix a dimension type $d$ in $\N I$, and 
fix representations $N$, $M$ in $\rd$. Let $E_N$ be a graded complementary 
of the tangent space of $\on$ at $N$ as in Proposition \ref{slice}. Then, the 
value of the polynomial
$\overline{\Omega_{M,N}}$ at $v^2$ equals the cardinality of the set 
$E_N\cap\om$.
\end{theorem}
\begin{proof} This is Proposition \ref{formula} combined with 
Corollary \ref{filtration}. 
\end{proof}

Note that the theorem implies the following curious identity :

\begin{corollary} Let $d$ be in $\N I$ and $N$ in $R_d$. Then,

$$\sum_P \overline{\Omega_{P,N}}=|E_N|=q^{\dim Ext^1(N,N)},$$
where $P$ runs over the set of isoclasses of representations of dimension type 
$d$.
\end{corollary}

We finish the section with the following remark. The (generalized) Hall 
polynomials are known to have leading coefficient equal to one. Hence, by the 
Lang-Weil theorem, all the varieties $E_N\cap\om$ have a unique irreducible 
component of maximal dimension.

\end{subsection}
\end{section}
\begin{section}{The preprojective variety and coefficients of the bar 
automorphism.}
\begin{subsection}{} 
For any arrow $\alpha 
\,:\,i\rightarrow j$ in the quiver $Q$, we define $i(\alpha)=i$ and 
$h(\alpha)=j$. Let $\qo$ be the opposite quiver, having the same vertices as 
$Q$, and an arrow $\alpha^*\,:\, j\rightarrow i$ for each arrow $\alpha 
\,:\,i\rightarrow j$ in $Q$.\par
Fix a dimension type $d=\sum_{i}d_ii$ in $\N I$. As above, $G_d$ acts on 
$\rd(Q)=\rd$ and 
on $\rd(\qo)$. Note that the map sending a linear map to its adjoint induces an 
isomorphim $\rd(Q)\rightarrow \rd(\qo)$, $M\mapsto M^*$.

Let $\Pi_d$ be the preprojective variety (see \cite{lusztig2}):
$$\Pi_d:=\{((M_\alpha)_\alpha,(N_{\alpha^*})_\alpha)\in\rd(Q)\times\rd(\qo),\,\mbox{ for all }i\in I:$$
$$\hskip 1 cm\sum_{\alpha\in Q_1, h(\alpha)=i}N_{\alpha^*}M_\alpha=
\sum_{\alpha\in Q_1, i(\alpha)=i}M_{\alpha}N_{\alpha^*}\}\subset
\rd(Q)\times \rd(\qo).$$
We denote by $p$ (resp.~$p^{op}$) the canonical projection $\Pi_d\rightarrow 
\rd(Q)$ (resp.~$\Pi_d\rightarrow \rd(\qo)$). We also define a twisted 
projection
$\tilde p$ : $\Pi_d\rightarrow \rd(Q)$ by mapping $(M,N)$ to $M+N^*$ using the 
identification $\rd(Q)\simeq\rd(\qo)$ above.
\end{subsection}
\begin{subsection}{}
Observe that, for all $N$ in $R_d$, the identification $\rd(Q)\rightarrow 
\rd(\qo)$
maps the orbit $G_d.N$ to $G_d.N^*$ and the tangent space $T_N(\on)$ to
$T_{N^*}({\mathcal O}_{N^*})$.\par
We consider the non degenerate pairing on $R_d(Q)\times R_d(\qo)$ given by
$$<(M_\alpha), (N_{\alpha^*})>=\sum_{\alpha\in Q_1}Tr(M_\alpha N_{\alpha^*}).$$
We have the following
\begin{lemma}
With respect to the pairing above, we have $T_N(\on)^{\bot}=\po(p^{-1}(N))$, for all $N$ in 
$R_d(Q)$.
\end{lemma}
\begin{proof}
By Lemma \ref{phi} (ii), an element $(f_{\alpha^*})$ of $R_d(\qo)$ belongs to 
$T_N(\on)^{\bot}$ if and only if for all $X$ in $\lgd$, we have
$$\sum_{\alpha\in Q_1}Tr((X_{h(\alpha)}N_\alpha-N_\alpha X_{i(\alpha)})
f_{\alpha^*})=0.$$
By well known properties of the trace form, this is equivalent to 
$$\sum_{\alpha\in Q_1}Tr(X_{h(\alpha)}N_\alpha f_{\alpha^*})- Tr (X_{i(\alpha)}
f_{\alpha^*}N_\alpha)=0,$$ thus

$$\sum_{i\in I}X_i(\sum_{i(\alpha)=i}N_\alpha f_{\alpha^*}-
\sum_{h(\alpha)=i}f_{\alpha^*}N_\alpha)=0.$$
Since the trace form is non-degenerate, this gives
$$\sum_{i(\alpha)=i}N_\alpha f_{\alpha^*}-
\sum_{h(\alpha)=i}f_{\alpha^*}N_\alpha=0,$$ which proves $(f_{\alpha^*})\in 
\po(p^{-1}(N))$ as required.
\end{proof}

For any morphism $\pi$ : $X\rightarrow Y$ of $k$-varieties, and any 
function $f$ (resp.~$g$) on $X$ (resp.~$Y$), we define as usual :
$$\pi^*(g)\,:\,X\rightarrow k,\; x\mapsto g(\pi(x)),$$
$$\pi_*(f)\,:\,Y\rightarrow k,\; y\mapsto \sum_{\pi(x)=y}f(x).$$
For a subset $A$ of $R_d$ or $R_d(\qo)$, we denote by $1_A$ the corresponding 
characteristic function. The previous lemma gives the following interpretation 
of the polynomials $\overline{\Omega_{M,N}}$ in terms of the geometry of the 
preprojective variety:
\begin{proposition}\label{tdual}
For $M$ in $R_d(Q)$ and $N^*$ in $R_d(\qo)$, we have $$(\po)_*(\tilde 
p)^*(1_{\om})(N^*)=\overline{\Omega_{M,N}}.$$
\end{proposition}
\begin{proof} Obviously, $(\po)_*(\tilde p)^*(1_{\om})$ belongs to ${\mathbb 
C}_{G_d}[R_d(\qo)]$.\par
Fix $N$ in $R_d$. By the definitions, $(\po)_*(\tilde p)^*(1_{\om})(N^*)$
can be rewritten as $(\po)_*(f_M)(N^*)$, where the function $f$ on $\Pi_d$
is defined by $f_M(A,B^*)=1$ if 
$A+B\in \om$, and 0 otherwise. Hence, by the previous lemma,
$$(\po)_*(\tilde p)^*(1_{\om})(N^*)=\sum_{A\in T_{N^*}({\mathcal 
O}_{N^*})^\bot}f_M(A,N^*)=|\om\cap(N+T_{N^*}({\mathcal O}_{N^*})^\bot)|.$$
By Theorem \ref{main}, this gives $(\po)_*(\tilde 
p)^*(1_{\om})(N^*)=\overline{\Omega_{M,N}}$.
\end{proof}
\end{subsection}
\end{section}

\end{document}